\documentclass[11pt]{article}

\input{amssym.def}
\input{amssym}
\usepackage{eufrak}

\newtheorem{theorem}{Theorem}

\newtheorem{corollary}[theorem]{Corollary}
\newtheorem{definition}[theorem]{Definition}

\newtheorem{remark}[theorem]{Remark}

\def\pa{{\partial}}

\def\un{\underline}
\def\qq{q^{-1}}

\def\Tr{\mathrm{Tr}}
\def\Trr{\Tr_R}

\def\LL{\hat{\cal L}}
\def\MM{{\cal{L}}}
\def\MM{{\cal{M}}}

\def\RR{R^{-1}}

\def\de{\delta}

\def\ot{\otimes}
\def\C{{\Bbb C}}

\def\ov{\overline}

\def\be{\begin{equation}}
\def\ee{\end{equation}}

\begin{document}

\makeatletter
\renewcommand{\theequation}{{\thesection}.{\arabic{equation}}}
\@addtoreset{equation}{section} \makeatother

\title{Matrix Capelli identities related to Reflection Equation algebra}
\author{\rule{0pt}{7mm} Dimitri
Gurevich\thanks{dimitri.gurevich@gmail.com}\\
{\small\it Universit\'e Polytechnique Hauts-de-France}\\
{\small\it F-59313 Valenciennes, France}\\
{\small \it and}\\
{\small \it Interdisciplinary Scientific Center J.-V.Poncelet}\\
{\small\it Moscow 119002, Russian Federation}\\
\rule{0pt}{7mm} Varvara Petrova\\
{\small\it
National Research University Higher School of Economics,}\\
{\small\it 20 Myasnitskaya Ulitsa, Moscow 101000, Russian Federation}\\
\rule{0pt}{7mm} Pavel Saponov\thanks{Pavel.Saponov@ihep.ru}\\
{\small\it
National Research University Higher School of Economics,}\\
{\small\it 20 Myasnitskaya Ulitsa, Moscow 101000, Russian Federation}\\
{\small \it and}\\
{\small \it
Institute for High Energy Physics, NRC "Kurchatov Institute"}\\
{\small \it Protvino 142281, Russian Federation}}

\maketitle

\begin{abstract}
By using the notion of Quantum Double we introduce analogs of partial derivatives on a Reflection Equation algebra, associated with a Hecke symmetry of $GL_N$ type.
We construct the matrix $L=MD$, where $M$ is the generating matrix of the Reflection Equation algebra and $D$ is the matrix composed of the quantum partial derivatives
and prove that the matrices $M$, $D$ and $L$ satisfy a matrix identity, called the matrix Capelli one.  Upon applying quantum trace, it becomes a scalar relation, which
is a far-reaching generalization of the classical Capelli identity.  Also, we get a generalization of the some higher Capelli identities defined in \cite{O} by A.Okounkov.
\end{abstract}

{\bf AMS Mathematics Subject Classification, 2010:} 81R50


\section{Introduction}

Let $M=\|m_i^j\|_{1\leq i,j \leq N}$ be a matrix with commutative entries and $D=\|\pa_i^j\|_{1\leq i,j \leq N}$ be the matrix composed of the  partial derivatives\footnote{Note that
$\pa_i^j{m^s_k}=\de_i^s\, \de_k^j$. Usually, in the Capelli identity one employs the matrix, transposed to our $D$.} $\pa_i^j=\partial/\partial m_j^i$. The famous Capelli identity reads
\be
\mathrm{cdet} (MD+K)=\det M\det D,
\label{cap}
\ee
where $\mathrm{cdet}$ is the so-called column-determinant and $K$ is a diagonal matrix of the form: $K=\mathrm{diag}(N-1, N-2,\dots,1, 0)$.

There are known many generalizations of this identity. We only mention the paper \cite{NUW}, where a quantum version of the Capelli identity was established, related to the
Quantum Group (QG) $U_q(sl_N)$ and its dual  algebra.

In the present note we exhibit another quantum version of the Capelli identity, which by contrast with \cite{NUW} is related to Reflection Equation (RE) algebras. By definition,
an RE algebra is a unital associative algebra $\MM(R)$ generated by entries of the matrix $M=\|m_i^j\|_{1\leq i,j \leq N}$
subject to the following  relation:
\be
R\,M_1R\,M_1-M_1R\,M_1R=0, \quad M_1 = M\otimes I,
\label{RE}
\ee
where $I$ is the unit  matrix and $R$ is a Hecke symmetry. The matrix $M$ is called the {\it generating matrix} of the algebra $\MM(R)$.

Let us precise that by a {\it Hecke symmetry} we mean a  braiding, meeting the Hecke condition:
$$
(q\, I\otimes I-R)( \qq I\otimes I+R)=0,\quad q\not\in\{0,\pm 1\},
$$
whereas by a {\it braiding} we mean  a solution of the braid relation:
$$
R_{12}\,R_{23}\,R_{12}=R_{23}\,R_{12}\,R_{23},\quad R_{12}=R\ot I,\quad R_{23}=I\ot R.
$$
Hereafter, $R$ is treated to be an $N^2\times N^2$ numerical matrix.

The best known examples of the Hecke symmetries are those coming from the QG $U_q(sl_N)$. These Hecke symmetries are deformations of the usual flips $P$.
Nevertheless, there exist other Hecke symmetries possessing this property (for instance, the Crammer-Gervais symmetries) as well as those which are not deformations of the usual
(or super-)flips.

We impose two additional requirements on the Hecke symmetry $R$: it should be skew-in\-ver\-ti\-ble and even (see \cite{GS1}). In such a case $R$ will be called the $GL_N$ type
Hecke symmetry. Note that if $R$ is a $GL_N$ type Hecke symmetry, then for generating matrix $M$ of the RE algebra $\MM(R)$ one can define the quantum (or $R$-)trace
$\Trr M$ and the quantum determinant ${\det}_R M$.

Besides,  for any $GL_N$ type Hecke symmetry $R$ we define analogs of the partial derivatives $\pa_i^j$ in such a way that the matrix $L=MD$, where
$D=\|\pa_i^j\|_{1\leq i,j \leq N}$, meets the relation:
\be
R\, L_1R\, L_1-L_1R\, L_1R=R\, L_1-L_1R.
\label{mRE}
\ee
An algebra $\LL(R)$, generated by entrees of the matrix $L=\|l_i^j\|_{1\leq i,j \leq N}$ is called a {\it modified RE algebra}. Note that as $R\to P$, the algebra $\MM(R)$
tends to $\mathrm{Sym}(gl_N)$, whereas the algebra $\LL(R)$ tends to $U(gl_N)$. This is one of the reasons why we consider the algebras $\MM(R)$
(resp., $\LL(R)$) for any $GL_N$ type Hecke symmetry $R$ as a quantum (or $q$-)analog of $\mathrm{Sym}(gl_N)$ (resp., $U(gl_N)$) .

Note that for the generating matrix $L$ of the algebra $\LL(R)$ the quantum trace $\Trr L$ and the quantum determinant ${\det}_R L$ are defined in the same way as for
the matrix $M$.  Namely, the quantum determinant of $L=M D$ with a proper shift  enters our quantum Capelli identity. It should be emphasized that this identity is
valid for the whole class of RE algebras $\MM(R)$, associated with $GL_N$ type Hecke symmetries $R$. Note that if $R\to P$ in the limit $q\to 1$ our quantum Capelli identity
turns into the classical one expressed as  in \cite{O}.

The note is organized as follows. In section 2 we exhibit the quantum double (QD) construction enabling us to introduce $q$-analogs of the partial derivatives in the entries of the
matrix $M$. In Theorem \ref{th:1} we present the matrix factorization identities which are called the {\it matrix Capelli identities}. Upon applying the $R$-trace, they turn into a quantum
version of the Capelli identity and some its generalizations which are the quantum counterparts of the {\em higher Capelli identities} (see Theorem in \cite{O}), corresponding to
one-column and one-row Young diagrams. In section 3 we give a proof of these identities. In section 4 we reduce the Capelli identity to a more conventional form, based on the
use of quantum determinants. Also, we compare our version of the Capelli identity with that from the article \cite{NUW}.

\section{Quantum partial derivatives and matrix Capelli identities}

In this section we deal with a skew-invertible Hecke symmetry  $R$ without assuming it to be even.

Consider two unital associative algebras $A$ and $B$ equipped with an invertible linear map $\sigma :A \ot B\to B\ot A$ which satisfies the following relations:
$$
\sigma\circ(\mu_{A}\otimes\mathrm{id}_B) =(\mathrm{id}_B\otimes\mu_{A})\circ\sigma_{12}\circ\sigma_{23} \quad{\rm on}\quad A\ot A\ot B,
$$
$$
\sigma\circ(\mathrm{id}_A\otimes\mu_{B}) =(\mu_{B}\otimes \mathrm{id}_A)\circ\sigma_{23}\circ\sigma_{12} \quad{\rm on}\quad A\ot B\ot B,
$$
$$
\sigma(1_A\ot b)=b\ot 1_A,\quad \sigma(a\ot 1_B)=1_B\ot a \qquad\forall\, a\in A,\, \forall\,b\in B,
$$
where $\mu_A: A\ot A\to A $ is the product in the algebra $A$, $1_A$ is its unit, and similarly for $B$. We call the data $(A,B, \sigma)$ a quantum double, if the map
$\sigma$ is defined in terms of a braiding $R$ (see \cite{GS2} for more detail).

Also, the map $\sigma$ defines permutation relations $a\ot b=\sigma (a\ot b)$, $a\in A$, $b\in B$ and due to this fact $\sigma$ is referred to as the {\it permutation map}.
If the algebra $A$ is equipped with a counit $\varepsilon:A\to \C$, then it becomes possible to define an action of the algebra $A$ onto $B$.

Below we deal with the QD $(A,B,\sigma)$,  where $B=\MM(R)$ with the generating matrix $M$ obeying (\ref{RE}), the algebra $A= {\cal D}(R^{-1})$ is the RE algebra
with the generating matrix $D = \|\partial_i^j\|$ satisfying the relation\footnote{The matrix $R^{-1}$ is also a Hecke symmetry but with $q$ replaced by $\qq$.}
\be
R^{-1}D_1R^{-1}D_1-D_1R^{-1}D_1R^{-1}=0
\label{DD}
\ee
and the permutation map is
$$
\sigma:\quad D_1R\, M_1R\to R\, M_1R^{-1}D_1+R \, 1_B1_A.
$$
Below we omit the factors $1_A$ and  $1_B$. The corresponding permutation relations can be written in the form:
\be
D_1R\, M_1 = R\, M_1 R^{-1} D_1 R^{-1}+ I.
\label{perr}
\ee
\begin{remark}
\rm
The quantum double $({\cal D}(R^{-1}), \MM(R), \sigma)$ with the permuttation relations (\ref{perr}) was obtained in \cite{GPS} from
the representation theory of the RE algebra.
\end{remark}
The permutation relations (\ref{perr}) are compatible with the associative structures of the both algebras $\MM(R)$ and
${\cal D}(R^{-1})$. To prove this we introduce the matrix notation:
$$
M_{\ov 1}=M_{\un 1}=M_1,\qquad M_{\ov {i+1}}=R_i M_{\ov i}R_i^{-1},\qquad M_{\un {i+1}}=R_i^{-1}M_{\un i}R_i,\quad i\ge 1,
$$
where $R_i:=R_{i\, i+1}:=I^{\otimes (i-1)}\otimes R\otimes I^{\otimes (p-i-1)}$ is an embedding of $R$ into the space of $N^p\times N^p$ matrices for
any $p\ge i+1$. Then the braid relation on $R$ allows one to prove the equivalence of two forms of defining relations of RE algebra $\MM(R)$:
$$
R\,M_1R\,M_1 - M_1R\,M_1R = 0\quad \Leftrightarrow\quad R\,M_{1}M_{\ov 2}-M_{1}M_{\ov 2}\,R =0.
$$
Note that the defining relations of RE algebra can be also written in terms of any higher copies of matrix $M$:
\be
R_p \,M_{\ov p}\,M_{\ov{p+1}}-M_{\ov p}\,M_{\ov{p+1}}\,R_p =0, \qquad \forall\,p\ge 1.
\label{h-copy}
\ee
By a straightforward calculation with the use of (\ref{perr}) we find:
\be
D_1\left(R_2M_{\ov 2}M_{\ov 3}-M_{\ov 2}M_{\ov 3}R_2 \right)=\left(R_2M_{\ov 2}M_{\ov 3}-M_{\ov 2}M_{\ov 3}R_2 \right) D_1\RR_1 R^{-2}_2\RR_1.
\label{ideal}
\ee
The relation (\ref{ideal}) entails that the defining ideal of the algebra $\MM(R)$ is preserved by the permutation relations. In a similar way it is possible to check that
the defining ideal of the algebra ${\cal D}(R^{-1})$ is also preserved by the permutation relations.

In order to get an action of the algebra ${\cal D}(R^{-1})$ onto $\MM(R)$ we introduce a counit in the algebra $A = {\cal D}(R^{-1})$ in the classical way:
$$
\varepsilon(1_A)=1_{\C},\quad \varepsilon(\pa_i^j)=0 \quad \forall\, i,j,\quad \varepsilon(a_1 a_2)=\varepsilon(a_1)\,\varepsilon(a_2) \quad \forall\, a_1, a_2 \in A.
$$
With this counit the action of $\partial_i^j$ on the generators $m_s^k$ reads:
$$
D_1\triangleright M_{\overline 2} = R_{12}^{-1}.
$$
The permutation relations (\ref{perr}) together with the counit map allow one to extend this action on the whole algebra $\MM(R)$. The elements $\partial_i^j$ with the
above action are treated as the quantum analogs of the usual partial derivatives in the commutative variables $m_i^j$. As was mentioned above, this action is
compatible with the algebraic structure of ${\cal M}(R)$ (see (\ref{ideal})). However, below we do not use the operator treatment of the quantum partial derivatives $\partial_i^j$.

\begin{remark}
\rm If $R$ is a Hecke symmetry coming from $U_q(sl_N)$ (the so-called Drin\-feld-Jim\-bo $R$-mat\-rix) then at the classical limit $q\rightarrow 1$
the permutation relations (\ref{perr}) turns into the usual Leinbniz rule for the commutative partial derivatives $\partial_i^j = \partial/\partial m_j^i$.
\end{remark}

With any Hecke symmetry $R$ we associate the idempotents $A^{(k)}$ and $S^{(k)}$ called the $R$-skew-symmetrizers and $R$-symmetrizers respectively. They are defined
by the following recursion:
\begin{eqnarray}
A^{(1)}=I,&\quad& A^{(k)}_{1\dots k}=\frac{1}{k_q}\, A^{(k-1)}_{1\dots k-1}\left(q^{(k-1)}\,I^{\otimes k}-(k-1)_q\, R_{k-1}\right) A^{(k-1)}_{1\dots k-1}, \quad  k\ge 2.
\label{skew}\\
S^{(1)}=I,&\quad& S^{(k)}_{1\dots k}=\frac{1}{k_q}\, S^{(k-1)}_{1\dots k-1}\left(q^{-(k-1)}\,I^{\otimes k}+(k-1)_q\, R_{k-1}\right) S^{(k-1)}_{1\dots k-1}, \quad  k\ge 2.
\nonumber
\end{eqnarray}
If $R$ is a $GL_N$ type Hecke symmetry, then $\mathrm{dim\,\, Im}\,A^{(N)} = 1$ and $A^{(N+1)}\equiv 0$.

Now we are ready to formulate the main result of the paper. We establish a series of matrix factorization identities which leads to the quantum versions
of the Capelli identity and some its generalizations called by A.Okounkov the ``higher Capelli identities'' in \cite{O}.

\begin{theorem} \label{th:1}
Let $L=MD$, where $M$ and $D$ are the generating matrices of the algebras $\MM(R)$ and ${\cal D}(\RR)$ from the quantum double defined by {\rm (\ref{RE})},
{\rm (\ref{DD})} and {\rm (\ref{perr})}. Then the following matrix factoriazation identities take place for $\forall\,k\ge 1$:
{\rm
\be
A^{(k)}L_{\overline 1}\,(L_{\ov 2}+q I)\dots (L_{\ov k}+q^{k-1}(k-1)_qI\,) \,A^{(k)}=
q^{k(k-1)}A^{(k)}M_{\overline 1}\dots M_{\ov k}\, D_{\ov k}\dots D_{\overline 1}
\label{th}
\ee}
{\rm
\be
S^{(k)} L_{\overline 1}\left(L_{\ov 2}-\frac{1}{q} I\right)\dots \left(L_{\ov k}-\frac{(k-1)_q}{q^{k-1}}I\right)S^{(k)} =
q^{-k(k-1)} S^{(k)}M_{\overline 1}\dots M_{\ov k}\, D_{\ov k}\dots D_{\overline 1}.
\label{th-s}
\ee}
\end{theorem}

Observe that Theorem \ref{th:1} is valid for {\it any} skew-invertible Hecke symmetry $R$. In the case when $R$ is a $GL_N$ type symmetry, the right hand side of
(\ref{th}) for $k=N$ can be presented as the product of quantum determinants of the matrices $M$ and $D$  (see the last section).

\begin{definition} Let $M$ be the generating matrix of an RE algebra $\MM(R)$. The quantities\footnote{Hereafter, we use the notation
 $ \langle X\rangle_{1\dots k}:= \mathrm{Tr}_{R(1\dots k)} X:=\mathrm{Tr}_{R(1)}\dots \mathrm{Tr}_{R( k)} X$,
where $X$ is an $N^k\times N^k$ matrix.}
$$
e_k(M)=\langle A^{(k)}\, M_{\overline 1}M_{\overline 2}\dots M_{\ov k}\rangle_{1\dots k}
$$
are called the elementary ($q$-)symmetric polynomials in the matrix $M$.
\end{definition}

By definition,  the quantum determinants of the matrices $M$ and $D$ are proportional to the highest elementary symmetric polynomials $e_N$ (similarly to the
classical matrix analysis):
\begin{equation}
{\det}_R M:=q^{N^2} \langle A^{(N)} M_{\overline 1}M_{\ov 2}\dots M_{\ov N}\rangle_{1\dots N},\quad
{\det}_{\RR}D:=q^{N^2}\langle A^{(N)}\, D_{\ov N}\dots D_{\overline 1}\rangle_{1\dots N}.
\label{r-det-n}
\end{equation}
The normalizing factior $q^{N^2}$ is introduced to simplify the formulae below. Note that in the definition of ${\det}_{\RR}D$ the inverse order of the matrix copies
$D_{\overline k}$ is used. This is motivated by the relations $(\ref{DD})$ imposed on $D$.

So, as a corollary of Theorem \ref{th:1}, we have the following version of the generalized quantum Capelli identities.

\begin{corollary} Under the  assumption of Theorem {\rm \ref{th:1}} the following identities hold for $\forall\,k\ge 1$:\rm
{\rm
\be
\langle A^{(k)}L_{\overline 1}\,(L_{\ov 2}+q I)\dots (L_{\ov k}+q^{k-1}(k-1)_qI\,)\rangle_{1\dots k} =
q^{k(k-1)}\langle A^{(k)}M_{\overline 1}\dots M_{\ov k}\, D_{\ov k}\dots D_{\overline 1}\rangle_{1\dots k}
\label{cap-as}
\ee}
{\rm
\be
\langle S^{(k)} L_{\overline 1}\left(L_{\ov 2}-\frac{1}{q} I\right)\dots \left(L_{\ov k}-\frac{(k-1)_q}{q^{k-1}}I\right)\rangle_{1\dots k} =
q^{-k(k-1)} \langle S^{(k)}M_{\overline 1}\dots M_{\ov k}\, D_{\ov k}\dots D_{\overline 1}\rangle_{1\dots k}.
\label{cap-s}
\ee}
\end{corollary}
Formulae (\ref{cap-as}) and (\ref{cap-s}) are generalizations  of the higher Capelli identities from \cite{O}, corresponding to
one-column and one-row Young diagrams respectively.
\begin{corollary}
\label{cor:6}
Under the  assumption of Theorem {\rm \ref{th:1}} the following quantum Capelli identity holds true:\rm
\be
\langle  A^{(N)}L_{\overline 1}\,(L_{\ov 2}+qI)\dots (L_{\ov N}+q^{N-1} (N-1)_qI) \rangle_{1\dots N} =
q^{-N} {\det}_RM \, {\det}_{\RR}D .
\label{cap1}
\ee
\end{corollary}

In the last section we consider the quantum determinants in more detail and complete the proof of this Capelli identity.

\section{Proof of Theorem \ref{th:1}}

We  only prove the identity (\ref{th}). The identity (\ref{th-s}) can be proven in  the same way.

Let us apply the induction in $k$. The base of induction for $k=1$ is tautological. We assume the identity (\ref{th}) to be true up to $k-1$ for some integer $k\ge 2$.
Consider the  matrix:
$$
F(\alpha) = A^{(k)} L_1(L_{\ov 2}+q\, I)(L_{\ov 3}+q^2\, 2_q\,I)\dots (L_{\ov {k-1}}+q^{k-2}\, (k-2)_q\, I)(L_{\ov {k}}+\alpha\, I) A^{(k)},
$$
where $\alpha$ is a numerical parameter.

Since $A^{(k-1)}$ is a polynomial in $R_i$ for $i\le k-2$, then $L_{\ov k}\, A^{(k-1)}=A^{(k-1)}\, L_{\ov k}$ as a consequence of the braid relation on $R$. Using this fact as
well as the identity $A^{(k)}=A^{(k)} A^{(k-1)}= A^{(k-1)}A^{(k)}$, we can rewrite $F(\alpha)$ in the form:
$$
F(\alpha)=A^{(k)}\un{A^{(k-1)} L_1(L_{\ov 2}+q\, I)\dots (L_{\ov {k-1}}+q^{k-2}(k-2)_q\, I)A^{(k-1)}}(L_{\ov {k}}+\alpha\, I)\,A^{(k)}.
$$
We transform the underlined expression in accordance with the induction hypothesis and get:
\be
F(\alpha)=q^{(k-1)(k-2)} A^{(k)}M_1\dots M_{\ov{k-1}}\,D_{\ov{k-1}}\dots D_1(L_{\ov k}+\alpha\, I)A^{(k)}.
\label{odin}
\ee
It remains to check that for  $\alpha=q^{k-1}(k-1)_q$ the expression $F(\alpha)$ turns into the right hand side of (\ref{th}). Expanding the brackets in (\ref{odin}) we obtain:
$$
q^{-(k-1)(k-2)}F(\alpha)=A^{(k)}M_1\dots M_{\ov{k-1}}D_{\ov{k-1}}\dots D_1 L_{\ov k}A^{(k)}+\alpha  A^{(k)} M_1\dots M_{\ov{k-1}}
D_{\ov{k-1}}\dots D_1A^{(k)}.
$$

Now, in the first summand we permute step by step all factors $D_{\ov i}$ with the element $L_{\ov k}$. Taking into account that
$L_{\overline k} = R_{k-1\rightarrow 1}L_1R^{-1}_{1\rightarrow k-1}$, where $ R_{1\to m}^{\pm}:=R_1^{\pm}\dots R_m^{\pm}$ (and similarly for $R^{\pm}_{m\rightarrow 1}$),
we find at the first step:
$$
D_1 L_{\ov k}=D_1R_{k-1\to 1} M_1 D_1\RR_{1\to k-1}=R_{k-1\to 2}D_1 R_1 M_1 D_1\RR_{1\to k-1}.
$$
In the last expression we replace the product $D_1R_1M_1$ with the use of (\ref{perr}):
$$
R_{k-1\to 2}\underline{D_1 R_1 M_1} D_1\RR_{1\to k-1} = R_{k-1\to 1}M_1\underline{R_1^{-1} D_1\RR_1D_1}\RR_{1\to k-1}+ D_1\, R_{k-1\to 2}\, \RR_{1\to k-1},
$$
then we change $R_1^{-1}D_1R_1^{-1}D_1$ for $D_1R_1^{-1}D_1R_1^{-1}$ according to (\ref{DD}) and finally get:
$$
D_1 L_{\ov k} = L_{\ov k}\, D_1R_{k-1\to 2}R_1^{-2}\RR_{2\to k-1}+D_1 R_{k-1\to 2}R_1^{-1}\RR_{2\to k-1}.
$$

So, the first summand in the above expression for $q^{-(k-1)(k-2)}F(\alpha)$ takes the form:
\begin{eqnarray*}
A^{(k)} M_1\!\!\!\!\!&\dots &\!\!\!\!\!M_{\ov{k-1}}D_{\ov{k-1}}\dots D_1L_{\ov k}A^{(k)}=\\
(-q)^2 \!\!\!\!\!&A^{(k)}&\!\!\!\!\! M_1\dots M_{\ov{k-1}} D_{\ov{k-1}}\dots D_2L_{\overline k}\,D_1A^{(k)}+(-q)\,A^{(k)}M_1\dots M_{\ov{k-1}}D_{\ov{k-1}}\dots D_1A^{(k)}.
\end{eqnarray*}
To get this expression, we ``evaluate'' the chains of $R$-matrices on the rightmost $R$-skew-sym\-met\-ri\-zer $A^{(k)}$ in accordance with the rules:
\be
A^{(k)}R_i^{\pm 1}= R_i^{\pm 1} A^{(k)}=-q^{\mp 1} A^{(k)}, \qquad 1\leq \forall \,  i \leq k-1. \label{consum} \ee

At the second step we permute $D_{\ov 2}$ and $L_{\ov k}$.  In the same way as above we find:
$$
D_{\ov 2} L_{\ov k}=L_{\ov k}D_{\ov 2}R_{k-1\to 3} R_2^{-2}\RR_{3\to k-1} +  D_{\ov 2}R_{k-1\to 3}R_2^{-1}\RR_{3\to k-1}.
$$
Note that all $R$-matrices in this formula commute with $D_1$ and therefore they can be moved to the right $A^{(k)}$ and converted to powers of $q$ according to (\ref{consum}).

By  induction in $p$ one can prove the general formula:
$$
D_{\ov p}L_{\ov k}=L_{\ov k}D_{\ov p} R_{k-1\to p+1}R_p^{-2} \RR_{p+1\to k-1} +  D_{\ov p}R_{k-1\to p+1}R_p^{-1}\RR_{p+1\to k-1}.
$$
Here also all terms $R^{\pm 1}_i$ $i\ge p$ commute with $D_{\ov{p-1}}\dots D_1$ and can be evaluated at the $R$-skew-symmetrizer $A^{(k)}$.

Finally, we get the following formula:
\begin{eqnarray*}
q^{-(k-1)(k-2)}F(\alpha) \!\!\!\!&=&\!\!\!\! q^{2(k-1)} A^{(k)} M_1\dots M_{\ov{k-1}}M_{\ov{k}}\, D_{\ov{k}} D_{\ov{k-1}}\dots D_1A^{(k)}\\
\rule{0pt}{5mm}&+&\!\!\!\!\! (\alpha -q-q^3-\dots -q^{2k-3}) A^{(k)}M_1\dots M_{\ov{k-1}} \,D_{\ov{k-1}}\dots D_{\ov 1} A^{(k)},
\end{eqnarray*}
where we substituted $L_{\overline k} = M_{\overline k}D_{\overline k}$.

At last, by setting $\alpha=q+q^3+\dots + q^{2k-3}=q^{k-1}(k-1)_q$ we kill the second term and get:
\be
F(q^{k-1}(k-1)_q)=q^{k(k-1)} A^{(k)}M_1\dots M_{\ov{k-1}}M_{\ov{k}} \, D_{\ov{k}} D_{\ov{k-1}}\dots D_1A^{(k)}.
\label{F-fin}
\ee
To complete the proof it remains to note that due to algebraic relations (\ref{h-copy}) the $R$-skew-sym\-met\-ri\-zer $A^{(k)}$ commute with the chain of $M$-matrices
$$
A^{(k)}M_1M_{\overline 2}\dots M_{\overline k} = M_1M_{\overline 2}\dots M_{\overline k} \,A^{(k)},
$$
and the same is true for the corresponding chain of $D$ matrices. Since $A^{(k)}A^{(k)}= A^{(k)}$, then in the right hand side of (\ref{F-fin}) one can leave only one element
$A^{(k)}$:
$$
A^{(k)}M_1\dots M_{\ov{k}} \, D_{\ov{k}} \dots D_1A^{(k)} \equiv A^{(k)}M_1\dots M_{\ov{k}} \, D_{\ov{k}} \dots D_1\equiv M_1\dots M_{\ov{k}} \,
D_{\ov{k}} \dots D_1\,A^{(k)}.
$$
This completes the inductive proof of (\ref{th}).

\section{Some aspects of quantum determinants}

It should be emphasized that the order $m$ of the highest non-trivial skew-symmetrizer $A^{(m)}$ can be different from $N$,  where $N^2\times N^2$ is the matrix
size of $R$.

\begin{definition}\rm
We say that a skew-invertible Hecke symmetry $R$ is  of rank $m$ if the $R$-skew-symmetrizers (\ref{skew}) satisfy the condition:\rm
$$
\mathrm{dim\,\, Im}\,A^{(m)}(R) = 1, \qquad A^{(m+1)}(R)\equiv 0.
$$
\end{definition}

Note that Corollary \ref{cor:6} remains valid, if in (\ref{cap1}) we replace $N$ by $m$ assuming the initial $GL_N$ type Hecke symmetry $R$ to be of rank $m$.

Since $A^{(m)}$ is an idempotent and  $\mathrm{dim\,\, Im}\,A^{(m)}=1$,  there exist two tensors
$|u\rangle= \|u_{i_1i_2\dots \,i_m}\|$ and $\langle v| = \|v^{i_1i_2\dots\,i_m}\|$
such that
$$
{A^{(m)}}_{i_1\dots\, i_m}^{\,\,\,j_1\dots\, j_m} = u_{i_1\dots\, i_m} v^{j_1\dots \,j_m}\qquad \mathrm{ and}\qquad \sum_{{i}}v^{i_1\dots \,i_m}u_{i_1\dots \,i_m} = 1.
$$
Using the above "bra" and "ket" notations, we can present these formulae as follows:
\be
A^{(m)}  = |u\rangle\langle v|  \qquad \mathrm{and}\qquad \langle v| u\rangle = 1.
\label{sh-not}
\ee

The quantum determinant is defined as in (\ref{r-det-n}) but with $N$ replaced by $m$:
$$
{\det}_R M\:=q^{m^2}\langle A^{(m)} M_1 M_{\ov 2}\dots M_{\ov m}\rangle_{1\dots m}.
$$
With the use of (\ref{sh-not}) we can prove the following {\it matrix} identity:
\begin{eqnarray}
A^{(m)}M_1\dots M_{\overline m} = A^{(m)}M_1\!\!\!\!&\dots&\!\!\!\!  M_{\overline m}A^{(m)} \nonumber\\
&&\!\!\!= |u\rangle\langle v | M_1\dots M_{\overline m}|u\rangle\langle v | 
 = A^{(m)}\langle v | M_1\dots M_{\overline m}|u\rangle.
\label{matr-id}
\end{eqnarray}
Upon calculating the $R$-trace over all spaces and taking into account that $\langle A^{(m)}\rangle_{1\dots m} = q^{-m^2}$ (see  \cite{GS1}),
we find that the quantum determinant is actually given by the {\it usual} trace of the form:
$$
{\det}_RM = \langle v| M_1M_{\overline 2}\dots M_{\overline m} |u\rangle := \mathrm{Tr}_{(1\dots m)}(A^{(m)}M_1M_{\overline 2}\dots M_{\overline m}).
$$
As for the quantum determinant of the matrix $D$ we have:
 $$
{\det}_{\RR}D:= q^{m^2}\langle A^{(m)} D_{\ov m}\dots D_{\overline 2}D_1\rangle_{1\dots m} = \langle v| D_{\ov m}\dots D_{\overline 2}D_1)|u \rangle.
$$
Note that this quantum determinant is defined with the same tensors $|u\rangle$ and $\langle v|$ though the matrix $D$ is subject to the RE with $R$ replaced by $\RR$. 
It can be explained by the fact that all skew-symmetrizers are invariant with respect to the replacement $R\to \RR$ and $q\to \qq$.

Consider now the identity (\ref{th}) for the Hecke symmetry of rank $m$. With the use of (\ref{matr-id}) the matrix structure of the right hand side of (\ref{th}) for $k=m$  can be 
transformed as follows:
$$
q^{m(m-1)}
A^{(m)} M_1\dots M_{\ov m}\,D_{\ov m}\dots D_1 = 
q^{m(m-1)}A^{(m)}\,{\det}_RM{\det}_{R^{-1}}D.
$$
Finally, by calculating the $R$-trace over all spaces of the both sides of (\ref{th}), we come to the desired form (\ref{cap1}) of the right hand side of the quantum Capelli identity:
$$
\langle  A^{(m)}L_{\overline 1}\,(L_{\ov 2}+qI)\dots (L_{\ov m}+q^{m-1} (m-1)_qI) \rangle_{1\dots m} = q^{-m}{\det}_RM{\det}_{R^{-1}}D.
$$
This completes the proof of Corollary \ref{cor:6}.

As follows from the results of \cite{GS3}, if a given Hecke symmetry $R$ is a deformation of the usual flip $P$, each of the determinants entering the right hand
side of the (\ref{cap1}) can be written as column-determinant or row-determinant. We do not know whether it is possible to do the same with the left hand side of (\ref{cap1}).
Also observe, that if an {\it involutive symmetry} $R$ (i.e. such that $R^2=I$) is a limit of a Hecke symmetry $R(q)$ as $q\to 1$, the corresponding Capelli identity can be
obtained from (\ref{cap1}) by setting $q=1$. Thus, it looks like the Capelli identity from \cite{O}, but the skew-symmetrizers and quantum determinants should be adapted
to $R=R(1)$.

At conclusion, we want to shortly compare our result for quantum Capelli identity with that of the paper \cite{NUW}. The authors of that paper deal with another quantum
version of the Capelli identity, related to the  QG $U_q(sl_N)$ and the corresponding RTT algebra. As for our results, we are working with quite different quantum algebra --- the
RE algebra and different quantum derivatives. Besides, we do not restrict ourselves with the $U_q(sl_N)$ $R$-matrix, our results are valid for the wide class
of RE algebras defined via {\it arbitrary} skew-invertible Hecke symmetries.


\begin{thebibliography}{NUW}

\bibitem[GS1]{GS1} Gurevich D, Saponov P. From Reflection Equation Algebra to Braided Yangians, Proceedings of the 1st International Conference on Mathematical Physics,
Grozny, Russia, 2016. Springer Proceedings in Math. and Statistics V.273 (2018).

\bibitem[GS2]{GS2} Gurevich D, Saponov P.  Doubles of Associative algebras and their Applications,  Phys. of Particles and Nuclei Letters 17 N5 (2020) 774--778

\bibitem[GS3]{GS3} Gurevich D, Saponov P. Determinants in Quantum Matrix algebras and integrable systems,
Theor. Math. Phys. 207 (2021), 261--276

\bibitem[GPS]{GPS} Gurevich D., Pyatov P., Saponov P. Braided Weyl algebras and differential calculus on $U(u(2))$, J. Geom. Phys.  62 (2012) 1175--1188.

\bibitem[NUW]{NUW} Noumi M., Umeda T., Wakayama M. A quantum analogue of the Capelli identity and elementary differential calculus  on $GL_q(n)$, Duke Math. J.
76 (1994) 567--594.

\bibitem[O]{O} Okounkov A. Quantum immanants and higher Capelli identities, Transformation Goups, 1 (1996) 99--126.

\end{thebibliography}
\end{document}